# CAN YOU HEAR THE FRACTAL DIMENSION OF A DRUM?*

WALTER ARRIGHETTI†
*Electronic Engineering Department, Università degli Studi di Roma "La Sapienza", via Eudossiana 18, Rome, 00184, Italy,* www.die.uniroma1.it/strutture/labcem/.

GIORGIO GEROSA
*Electronic Engineering Department, Università degli Studi di Roma "La Sapienza", via Eudossiana 18, Rome, 00184, ITALY,* www.die.uniroma1.it/strutture/labcem/.

Electromagnetics and Acoustics on a bounded domain is governed by the Helmholtz's equation; when such a domain is [pre-]fractal described by means of a 'just-touching' Iterated Function System (IFS) spectral decomposition of the Helmholtz's operator is self-similar as well. Renormalization of the Green's function proves this feature and isolates a subclass of eigenmodes, called *diaperiodic*, whose waveforms and eigenvalues can be recursively computed applying the IFS to the initiator's eigenspaces. The definition of *spectral dimension* is given and proven to depend on diaperiodic modes only for a wide class of IFSs. Finally, asymptotic equivalence between box-counting and spectral dimensions in the fractal limit is proven. As the "self-similar" spectrum of the fractal is enough to compute box-counting dimension, positive answer is given to title question.

## 1. Introduction

### 1.1. *Between spectral and fractal geometry*

Marc Kac wondered in 1966 (in a famous paper entitled '*Can you hear the shape of a drum?*', [1]) whether the *shape* of a plane domain could be inferred from the sole spectrum of its Laplace's operator with Neumann/Dirichlet's boundary conditions, i.e. whether the shape of a membrane can be inferred by just "hearing" all its vibrating modes. This conjecture was confuted in 1992, when two plane domains were found to have the same spectrum, but different shapes [2]. Thereinafter Euclidean domains are usually split into equivalence classes of *isospectrality*, whose members have the same laplacian spectrum but different shapes. Computationally speaking, Laplace's operator on a bounded domain has a discrete spectrum (formed by at most a countable infinity of eigenvalues); the shape of a drum has, in general, a non-countable infinity of degrees of freedom.

† Email addresses: arrighetti@die.uniroma1.it and gerosa@die.uniroma1.it.

This problem can be extended of course to arbitrary, finite-dimensional domains and to self-similar/fractal sets. Title question aims to the extraction of a single number (the "fractal dimension") from that countable spectrum, that is computing a measure of the domain's complexity and self-symmetry.

### 1.2. *Brief review on Iterated Function Systems*

Let $\mathbf{w}_1,\mathbf{w}_2,\ldots,\mathbf{w}_p:\mathbb{R}^d\to\mathbb{R}^d$ be $p$ contractions with contraction ratios $c_j\in]0,1[$, $1\le j\le p$ respectively, i.e.:

$$\left\|\mathbf{w}_j(\mathbf{x}) - \mathbf{w}_j(\mathbf{y})\right\| \le c_j \left\|\mathbf{x} - \mathbf{y}\right\|, \qquad \forall \mathbf{x},\mathbf{y}\in\mathbb{R}^d,\ 1\le j\le p.$$

Set $\{\mathbf{w}_1,\mathbf{w}_2,\ldots,\mathbf{w}_p\}$ is said to be an *Iterated Function System* (*IFS*), or Iterated Similarity System (*ISS*) if all the $p$ contractions are similarities. It can be shown (cfr. [3], [6]) that, letting $\wp_{©}(S)$ be the set of all compact subsets of a metric space $S$, the contraction mapping $\mathbf{w}:\wp_{©}(\mathbb{R}^d)\to\wp_{©}(\mathbb{R}^d)$ defined as

$$C \mapsto \mathbf{w}(C) := \bigcup_{j=1}^{p} \mathbf{w}_j(C)\,, \quad \forall C\in\wp_{©}(\mathbb{R}^d) \tag{1}$$

admits one compact set $F\subset\mathbb{R}^d$ which is a fixed point $F\in\wp_{©}(\mathbb{R}^d)$ (respect to the Hutchinson's metric); $F$ is called the IFS' *attractor*. Thus for *any* compacts $E_0\in\wp_{©}(\mathbb{R}^d)$ and given the sets' sequence $(E_N)_{N\in\mathbb{N}}$ such that $\forall N\in\mathbb{N}_0$ $E_{N+1}=\mathbf{w}(E_N)=\mathbf{w}^N(E_0)$, the following limit exists:

$$F = \lim_N E_N \equiv \lim_N \mathbf{w}^N(E_0) = \bigcap_{N=0}^{\infty} E_N\,. \tag{2}$$

Attractor $F$ is often *self-similar* (i.e. it is made up of similar copies of itself at infinitely many smaller length-scales), does not depend on initiator $E_0$ and has non-integer fractal dimension*s* (e.g. Hausdorff and box-counting, [5]). Sets of sequence $(E_N)_{N\in\mathbb{N}}$ depend on initiator $E_0$ and might show self-similarity up to a finite length-scale only. $E_N$ is thus called $N^{th}$ *prefractal* of the IFS (of $E_0$ initiator).

An IFS is said to be disconnected if the $p$ copies are disjoint, or *'just-touching'* if they are not overlapping (but their boundaries may partially coincide), i.e. $\mathbf{w}_i(E)^\circ\cap\mathbf{w}_j(E)^\circ=\varnothing$, $1\le i<j\le p$.

### 1.3. *Box-counting dimension*

As previously said the main indicator of "fractality" is the concept of fractal dimensions, which several definitions exist for and whose values often coincide for IFSs' attractors. Just one will be considered then and will be used in §3.2.

Let $(S,\mathrm{d})$ be a metric structure, $\delta>0$ and $\{U_i\}_{i\in I}$ be a *$\delta$-covering* of $S$ (i.e. $\delta = \sup_{i\in I} \operatorname{diam} U_i$ and $\bigcup_{i\in I} U_i = S$); let $\mathbb{N}_\delta(S)$ be the minimal numer of sets $\delta$-covering $S$, no matter what kind of sets they are. The *box-counting dimension* of $S$, whenever it exists (i.e. whether the following limit exists), is defined as:

$$\dim_{\mathrm{B}} S := -\lim_{\delta\to 0} \frac{\log \mathbb{N}_\delta(S)}{\log \delta}\,. \tag{3}$$

Box-counting dimension does not depend from specific $\delta$-covering types (cfr. [5]), so its sets can belong to any subclasses, e.g. either overlapping balls and polytopes, pavements of (i.e. 'just-touching') polytopes, etc. A useful method to approximate it is the *coarse graining*, a hierarchy refinement of such coverings (usually packings of hypercubes) with a decreasing sequence of discrete diameters $(\delta_n)_{n\in\mathbb{N}}$, stopped whenever desired approximation is reached.

In the case of a 'just-touching' ISS, all the fractal dimensions coincide, depend on just the contraction ratios (cfr. §1.2) and solve the equation:

$$\sum_{j=1}^{p} c_j^{\dim_{\mathrm{B}} F} = 1\,. \tag{4}$$

If the similarities all have the contraction ratio $c$, then $\dim_{\mathrm{B}} S = -\log_c p$ .

## 2. Self-similar spectral decomposition of the Green's function

### 2.1. *Spectral decomposition for 'just-touching' prefractals*

Let $E_0\in\wp_{\copyright}(\mathbb{R}^d)$, let $E_N$ be its $N^{\text{th}}$ prefractal under the IFS of §1.2 and $G_N(\mathbf{x},\mathbf{x}';\lambda)$ be the $H_0^2(E_N)$ Green's function for the Helmholtz's operator $\nabla^2-\lambda\mathrm{id}$ on $E_N$ with zero boundary condition on $\partial E_N$ i.e., weakly ($\nabla^2$ operates on $\mathbf{x}\in E_N$),

$$\big\langle (\nabla^2 - \lambda\mathrm{id}) G_N \,\big|\, f \big\rangle_{H^2(E_N)} = f(\mathbf{x}')\,, \quad \forall f \in H_0^2(E_N)\,. \tag{5}$$

Being $E_N$ compact $\forall N\in\mathbb{N}_0$, Laplace's operator $\nabla^2$ is self-adjoint, non-positive definite, and compact; its spectrum $\operatorname{spec}\nabla^2=\{\lambda_{N,n}\}_{n\in\mathbb{N}}$ is countable, real, non-positive and $\lim_N \lambda_{N,n} = -\infty$, $\forall n\in\mathbb{N}$. Its eigenspaces are finite-dimensional and provide a complete basis for the Sobolev-Hilbert space $H_0^2(E_N)$ of $\mathrm{L}^2$ functions, zero on $\partial E_N$, with weak II-order derivatives:

$$H_0^2(E_N) = \bigoplus_{n\in\mathbb{N}} \mathrm{Ker}\left(\nabla^2 - \lambda_{N,n}\mathrm{id}\right) = \overline{\mathrm{span}\left\{\varphi_{N,n}\right\}_{n\in\mathbb{N}}} \,. \tag{6}$$

Eigenfunction associated to eigenvalue $\lambda_{N,n}$ (counting multiplicities) is $\varphi_{N,n}$.

As the IFS is 'just-touching' and domain $E_N$=$\mathbf{w}(E_{N-1})$ is the union of $\mathbf{w}_1(E_{N-1})$, $\mathbf{w}_2(E_{N-1})$, …, $\mathbf{w}_p(E_{N-1})$, thanks to operator's linearity, all the eigenfunctions $\varphi_{N-1,n}$ associated to $(N-1)^{\text{th}}$ prefractal, once "contracted and copied" over the $p$ copies via $\mathbf{w}_j$ contractions, $1\leq j\leq p$, still are the $N^{\text{th}}$ prefractal's eigenfunctions, despite not forming a *complete basis*. These eigenfunctions are self-similar, i.e. they are similar copies of initiatior set $E_0$'s eigenfunctions, rescaled down to the $N^{\text{th}}$-order prefractal. They can be directly computed from eigenfunctions of the initiator and "periodically" copied on every rescaled copies of it "towards" the prefractal's $N^{\text{th}}$ order. That is why they are called *diaperiodic mode*s, [7]. Their zero-boundary condition is copied throughout the $p$ copies so diaperiodic eigenfunctions are zero on $\partial\mathbf{w}_1(E_{N-1})$, $\partial\mathbf{w}_2(E_{N-1})$, …, $\partial\mathbf{w}_p(E_{N-1})$, i.e. on the intersections of the previous prefractal's copies too. Unless the IFS is disconnected, additional modes will be present, nonzero on copies' intersections and referred to as *interconnective mode*s.

Scalar eigenfunctions associated to 0 eigenvalue are *finite*, different from $N$ to $N$ and usually depend on domain's *cohomology* i.e., roughly, on the number of its "holes". $0\in\mathrm{spec}\nabla^2$ whenever domains $E_N$ are simply connected. As domain's topology will not further investigated here [6], their contribute to the Green's function (just indicated as $\gamma_N/\lambda$, with a polar singularity in $\lambda$=0, cfr. §2.2), will not be treated.

A 3-plet of indices $(N,j,n)$ is used for both eigenfunctions and (nonzero) eigenvalues: $N$ is the prefractal order, $n$ is the usual index counting independent eigenfunctions (or eigenvalues, although they can be equal due to degeneracy), while $0\leq j\leq p$ separates the eigenfunctions such that $j$=0 refers to interconnective modes and $1\leq j\leq p$ refers to diaperiodic modes obtained rescaling the whole set of nonzero eigenmodes of the $(N-1)^{\text{th}}$ prefractal, copied to the $j^{\text{th}}$ copy $\mathbf{w}_j(E_{N-1})\subset E_N$ only and letting it to be zero everywhere else on $E_N$. Thus generic eigenfunction $\varphi_{N,j,n}$ , $1\leq j\leq p$, is:

$$\varphi_{N,j,n}(\mathbf{x}) := \varphi_{N-1,n}\left(\mathbf{w}_j(\mathbf{x})\right)\chi_{\mathbf{w}_j(E_{N-1})}(\mathbf{x}) \tag{7}$$

(characteristic functions $\chi_A$'s will be further suppressed for simplicity's sake).

Eigenfunctions' families for cohomologic and interconnective modes should be computed for each $N$, and might show self-similarity features too.

**2.2.** ***Renormalization of the Green's function***

From the considerations of §2.1, the Green's function $G_N$ can be decomposed as:

$$G_N(\mathbf{x},\mathbf{x}';\lambda) = \frac{\gamma_N(\mathbf{x},\mathbf{x}')}{\lambda} + g_N(\mathbf{x},\mathbf{x}';\lambda) = \\ = \frac{\gamma_N(\mathbf{x},\mathbf{x}')}{\lambda} + g_N^0(\mathbf{x},\mathbf{x}';\lambda) + g_N^1(\mathbf{x},\mathbf{x}';\lambda) \ , \tag{8}$$

where:

$$g_N(\mathbf{x},\mathbf{x}';\lambda) = \sum_{n\in\mathbb{N}} \frac{\varphi_{N,n}(\mathbf{x})\varphi_{N,n}(\mathbf{x}')}{\lambda - \lambda_{N,n}} \equiv \sum_{n\in\mathbb{N}} \sum_{j=0}^{p} \frac{\varphi_{N,j,n}(\mathbf{x})\varphi_{N,j,n}(\mathbf{x}')}{\lambda - \lambda_{N,j,n}} .$$

Function $g_N$ is the Green's function expunged from cohomologic modes (not present at all whenever $E_N$ is simply linearly connected), which is further decomposed in a "diaperiodic part" $g_N^0$ and an "interconnective part" $g_N^1$ :

$$g_N^0(\mathbf{x},\mathbf{x}';\lambda) = \sum_{n\in\mathbb{N}} \sum_{j=1}^{p} \frac{\varphi_{N,j,n}(\mathbf{x})\varphi_{N,j,n}(\mathbf{x}')}{\lambda - \lambda_{N,j,n}} ;$$

$$g_N^1(\mathbf{x},\mathbf{x}';\lambda) = \sum_{n\in\mathbb{N}} \frac{\varphi_{N,0,n}(\mathbf{x})\varphi_{N,0,n}(\mathbf{x}')}{\lambda - \lambda_{N,0,n}} .$$

Whenever an ISS is considered, eigenvalues get rescaled, as $N$ increases, as $\lambda_{N,j,n} \geq c_j^{-1}\lambda_{N-1,n}$ and this becomes an equality whenever $\mathbf{w}_1$, $\mathbf{w}_2$, …, $\mathbf{w}_p$ have homogeneous ratios $c_j$'s (i.e. they are combinations of homotheties, reflections, translations and rotations only). In the latter case $g_N^0$ is recursively renormalized [4] from the Green's function initiator $g_0^0 \equiv g_0$ and using (7):

$$\begin{aligned} g_N^0(\mathbf{x},\mathbf{x}';\lambda) &= \sum_{j=1}^{p}\sum_{n=1}^{\infty} \frac{\varphi_{N-1,n}\left(\mathbf{w}_j(\mathbf{x})\right)\varphi_{N-1,n}\left(\mathbf{w}_j(\mathbf{x}')\right)}{\lambda - c_j^{-1}\lambda_{N-1,n}} \\ &= \sum_{j=1}^{p} c_j \sum_{n=1}^{\infty} \frac{\varphi_{N-1,n}\left(\mathbf{w}_j(\mathbf{x})\right)\varphi_{N-1,n}\left(\mathbf{w}_j(\mathbf{x}')\right)}{c_j\lambda - \lambda_{N-1,n}} \\ &= \sum_{j=1}^{p} c_j g_{N-1}\left(\mathbf{w}_j(\mathbf{x}),\mathbf{w}_j(\mathbf{x}');c_j\lambda\right) . \end{aligned} \tag{9}$$

Figure 1 shows a sample of eigenmodes for a *Šerpinskij carpet*'s 1st-step prefractal. Left to right: combination of a first 'vertical' diaperiodic mode and a second 'horizontal' interconnective mode; first square diaperiodic mode; second square diaperiodic mode. Rightmost image shows the *first* square diaperiodic mode for the *2nd-step* prefractal. Graphs over each image represent the wavelengths' distribution along a vertical section of the prefractals.

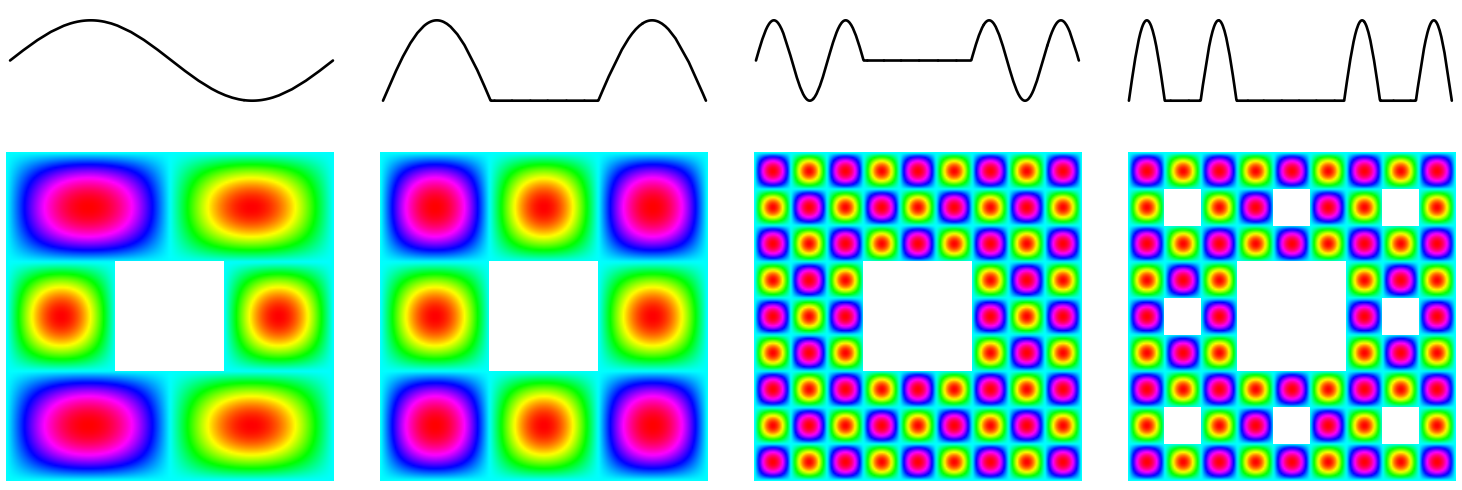

Figure 1. Samples of diaperiodic modes (and their symbolic wavelengths) on a square Šerpinskij carpet's 1st- and 2nd-step prefractal.

## 2.3. *Self-similar eigenvalues' scaling*

Renormalization (9) also shows the scaling of eigenvalues' distribution as the IFS ordering increases: at $N^{\text{th}}$ iterate the $(N-1)^{\text{th}}$-step's whole spectrum is rescaled $p$ times by coefficient $c_j$; then $N^{\text{th}}$-step interconnective eigenvalues are added. Figure 2 shows an example of a *Šerpinskij gasket* prefractal's self-similar spectrum.

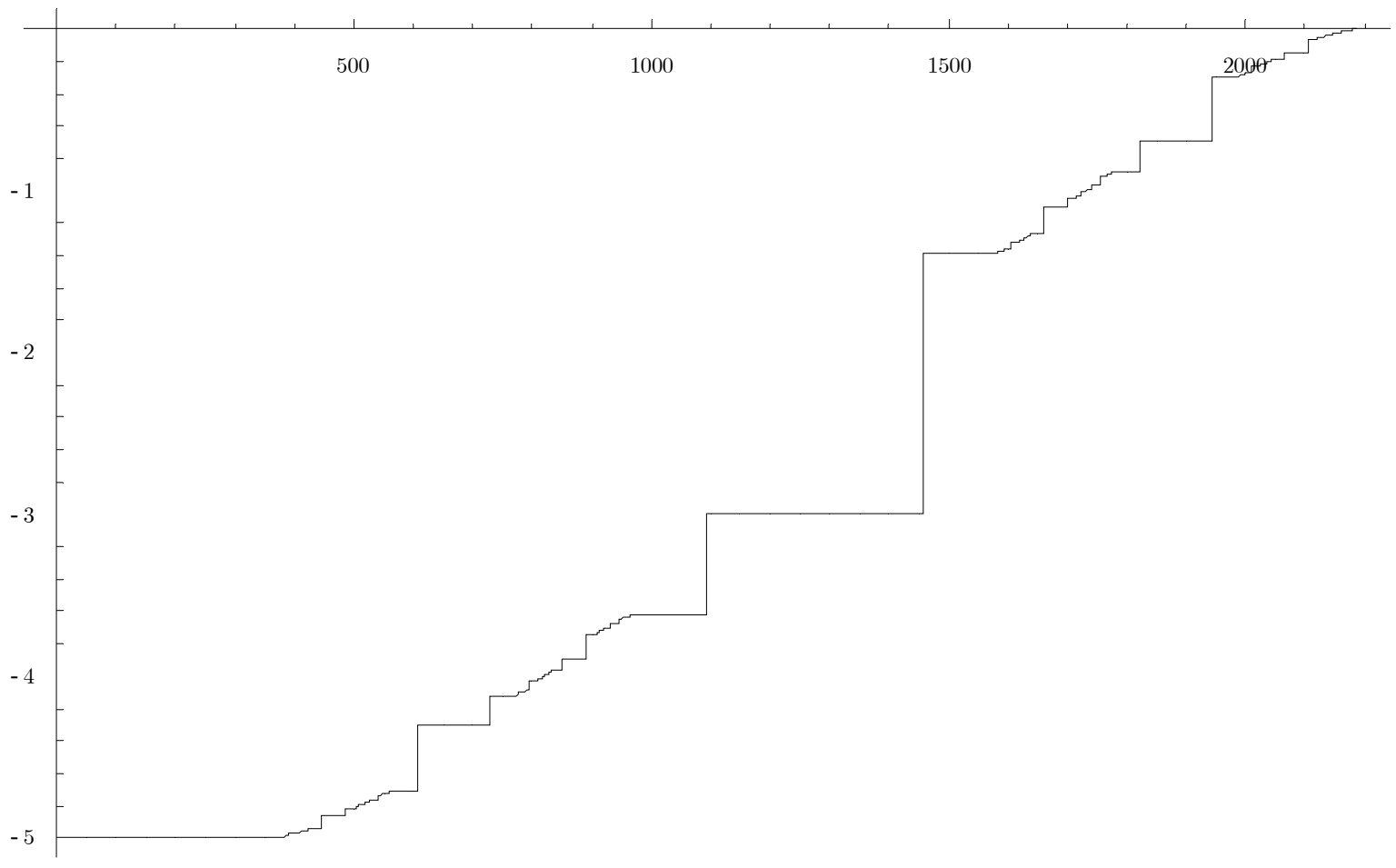

Figure 2. Laplacian eigenvalues' distribution (together with their multiplicities) of a Šerpinskij gasket's 7th-step prefractal [6].

Eigenvalues' count is reported in abscissa, versus eigenvalues; each "*plateau*" in the graph stands for a multiple eigenvalue (the boarder the plateau, the higher the eigenvalue's multiplicity). Boarder plateaux are associated to bigger triangular holes in the gasket's prefractal, so the self-similar geometry of the set is reflected in the self-similar distribution of its eigenvalues.

## 3. One can hear the fractal dimension of a drum

### 3.1. *Spectral dimension and its asymptotics*

Let $E \subset \mathbb{R}^d$ be a compact $d$-dimensional domain and $\nabla^2$ be its (scalar) laplacian operator with zero boundary conditions. The *spectral dimension* of set $E$ is defined as:

$$\dim_\Delta E := \frac{\displaystyle\sum_{\lambda \in \operatorname{spec} \nabla^2 \setminus \{0\}} \log \operatorname{mul} \lambda}{\displaystyle\sum_{\lambda \in \operatorname{spec} \nabla^2 \setminus \{0\}} \log |\lambda|} . \tag{10}$$

It exists for every 'just-touching' ISSs' prefractal and their attractor $F$ with $\dim_\Delta F$ computed as the limit of its prefractals' spectral dimension as $N \to \infty$:

$$\dim_\Delta F = \lim_N \dim_\Delta E_N .$$

Let $1 < q \leq p$ be the number of the ISS' similarities with different contraction ratios, $c_1$, $c_2$, …, $c_q$ (case $q=1$ is 'easier', slightly different and will be treated later); so there are $p_1$ similarities with contraction ratio $c_1$, $p_2$ similarities with contraction ratio $c_2$, …, $p_q$ similarities with contraction ratio $c_q$; trivially $p_1+p_2+\ldots+p_q=p$. In this case $\operatorname{mul}\lambda_{N,j,n} = p_j \operatorname{mul}\lambda_{N-1,n}$ , so the numerator of (10) can be rescaled as:

$$\begin{aligned}
\log \mathrm{mul}\lambda_{N,n} &= \sum_{j=0}^{q} \log \mathrm{mul}\lambda_{N,j,n} = \\
&= \log \mathrm{mul}\lambda_{N,0,n} + \sum_{j=1}^{q} \log\left( p_j \mathrm{mul}\lambda_{N-1,n} \right) \\
&= \log \mathrm{mul}\lambda_{N,0,n} + q\sum_{j=0}^{q} \log \mathrm{mul}\lambda_{N-1,j,n} + \sum_{j=1}^{q} \log p_j \\
&= \log \mathrm{mul}\lambda_{N,0,n} + q \log \mathrm{mul}\lambda_{N-1,0,n} + \\
&\quad + q\sum_{j=1}^{q} \log\left( p_j \mathrm{mul}\lambda_{N-2,n} \right) + \sum_{j=1}^{q} \log p_j \\
&= \log \mathrm{mul}\lambda_{N,0,n} + q \log \mathrm{mul}\lambda_{N-1,0,n} + \\
&\quad + q^2 \log \mathrm{mul}\lambda_{N-2,n} + (1+q)\sum_{j=1}^{q} \log p_j = \dots \\
\dots &= q^N \log \mathrm{mul}\lambda_{0,n} + \sum_{k=0}^{N-1} q^k \left( \log \mathrm{mul}\lambda_{N-k,0,n} + \sum_{j=1}^{q} \log p_j \right);
\end{aligned} \tag{11}$$

denominator of (10) is rescaled as:

$$\begin{aligned}
\sum_{j=0}^{p} \log \lambda_{N,j,n} &= \log \lambda_{N,0,n} + \sum_{j=1}^{q} \log \frac{\lambda_{N-1,n}}{c_j} \\
&= \log \lambda_{N,0,n} + q\sum_{j=0}^{q} \log \lambda_{N-1,j,n} - \sum_{j=1}^{q} \log c_j \\
&= \log \lambda_{N,0,n} + q \log \lambda_{N-1,0,n} + q^2 \log \lambda_{N-2,0,n} + \\
&\quad - (1+q)\sum_{j=1}^{q} \log c_j = \dots \\
\dots &= q^N \log \lambda_{0,n} + \sum_{k=1}^{N} \left( q^{N-k} \log \lambda_{k,p(p+1)+n} - q^{k-1} \cdot \sum_{j=1}^{q} \log c_j \right).
\end{aligned} \tag{12}$$

Putting together (11) and (12), one has:

$$\dim_{\Delta} E_N = \frac{\sum_{n=1}^{\infty} \log \mathrm{mul}\, \lambda_{N,n}}{\sum_{n=1}^{\infty} \log \left| \lambda_{N,n} \right|} = \tag{13}$$

$$= \frac{\sum_{k=0}^{N-1} q^k \cdot \sum_{j=1}^{q} \log p_j \cdot \sum^{\infty} 1 + \sum_{n=1}^{\infty} \left( q^N \log \mathrm{mul} \lambda_{0,n} + \sum_{k=0}^{N-1} q^k \log \mathrm{mul} \lambda_{N-k,0,n} \right)}{\sum_{k=0}^{N-1} q^k \cdot \sum_{j=1}^{q} \log \frac{1}{c_j} \cdot \sum^{\infty} 1 + \sum_{n=1}^{\infty} \left( q^N \log \lambda_{0,n} + \sum_{k=1}^{N-1} q^{N-k} \log \lambda_{k,p(p+1)+n} \right)}$$

Rightmost sums at numerator and denominator, diverging as $\mathrm{O}(q^{N-1})$, can be neglected respect to the other sums, increasing as $\mathrm{O}(q^N)$, thus proving that:

$$\dim_\Delta F = \lim_N \dim_\Delta E_N \sim$$

$$\sim \lim_N \frac{\frac{q^N - 1}{q-1} \sum_{j=1}^{q} \log p_j \cdot \sum^{\infty} 1 + q^N \sum_{n=1}^{\infty} \log \mathrm{mul} \lambda_{0,n}}{\frac{q^N - 1}{q-1} \sum_{j=1}^{q} \log \frac{1}{c_j} \cdot \sum^{\infty} 1 + q^N \sum_{n=1}^{\infty} \log \lambda_{0,n}} < +\infty\,.$$

This result also proves that the *spectral dimension of a* [*pre-*]*fractal only depends on* [*diaperiodic eigenvalues and*] *the ISS' contraction ratios*.

Case $q$=1 ($p$ similarities with one ratio $c\in]0,1[$) is not only proven to be well-posed, but it also easily equals spectral dimension with box-counting dimension (cfr. §1.3). Estimating asymptotics as (13) leads to:

$$\dim_\Delta F \sim \lim_N \frac{\frac{N(N+1)}{2} \log p \cdot \cancel{\sum^{\infty} 1} + \cancel{q^N \sum_{n=1}^{\infty} \log \mathrm{mul} \lambda_{0,n}}}{\frac{N(N+1)}{2} \log \frac{1}{c} \cdot \cancel{\sum^{\infty} 1} + \cancel{q^N \sum_{n=1}^{\infty} \log \lambda_{0,n}}} \tag{14}$$

$$= \lim_N \frac{\log p}{\log c^{-1}} \equiv -\frac{\log p}{\log c} \equiv \dim_{\mathrm{B}} F.$$

Case $q$>1 is left with expression (13), which has no general closed-form sums, so it is better to check whether it equals some other known quantities.

### 3.2. *Spectral and box-counting dimensions*

Due to the contractive nature of an ISS, an initiator $B_0 \subset \mathbb{R}^d$ is chosen such that its prefractals' sequence $(B_N)_{N\in\mathbb{N}} \to F$ can be used as a coarse graining for computing $\dim_{\mathrm{B}} F$. $B_0$ can be chosen as a suitable hypercube [packing of hypercubes] of [maximal] side $\delta_0$; prefractals $B_N$'s, via IFS iteration, are made of hypercubes of maximal side $\delta_N$ (of decreasing diameters $\delta_N \sqrt{d} \to 0$ ), yet 'just-touching' with each other.

The laplacian spectrum of such hypercubes is well known to be the set:

$$\left\{ -\frac{\pi \,\|\mathbf{m}\|}{\delta_N} \;\middle|\; \forall \mathbf{m} \in \mathbb{N}^d \right\}; \qquad \mathrm{mul}\left( -\frac{\pi \,\|\mathbf{m}\|}{\delta_N} \right) \sim \mathrm{O}(d!)\,.$$

As $\dim_\Delta F$ does not *asymptotically* dependent on interconnecrtive eigenvalues (cfr. §3.1), just the algebraic multiplicities count for $\mathbb{N}_{\delta_N}(B_N)$. Considering the coarse graining properties (cfr. §1.3 and [5]), one gets:

$$\begin{aligned}
\dim_\Delta F = \lim_N \dim_\Delta E_N &= \lim_N \frac{\displaystyle\sum_{\|\mathbf{m}\| \in \mathbb{N}_0^d \setminus \{\mathbf{0}\}} \log \mathrm{mul}\, \frac{\pi \,\|\mathbf{m}\|}{\delta_N}}{\displaystyle\sum_{\|\mathbf{m}\| \in \mathbb{N}_0^d \setminus \{\mathbf{0}\}} \log \frac{\pi \,\|\mathbf{m}\|}{\delta_N}} \\
&\sim \lim_N \frac{\log \mathbb{N}_{\delta_N}(B_N) + \cancel{\log d!}}{\cancel{\log \pi} - \log \delta_N \cancel{\sqrt{d}}} = -\lim_N \frac{\log \mathbb{N}_{\delta_N}(F)}{\log \delta_N} \qquad (15) \\
&= -\lim_{\delta \to 0} \frac{\log \mathbb{N}_\delta(F)}{\log \delta} \equiv \dim_{\mathrm{B}} F\,.
\end{aligned}$$

So the spectral dimension and the box-counting dimension of such an ISS' attractor coincide and this gives the «*yes, one can*» answer to title question.

### 3.3. *Examples*

Simple examples of spectral dimension's computation, for both standard sets and ISS' attractors are given.

[0,1] $\subset \mathbb{R}$ is considered first. The unit interval is the trivial attractor to a 'just-touching' ISS made up of $w_1(x) = \frac{1}{2}x$ and $w_2(x) = \frac{1}{2}(x+1)$. Associated prefractals all trivially equal [0,1] itself and just the $N^{\text{th}}$-step nonzero diaperiodic eigenvalues $\{-2^N n\pi\}_{n\in\mathbb{N}}$ count, with $\mathrm{mul}(2^N n\pi)=2^N$:

$$\begin{aligned}
\dim_\Delta[0,1] &= \lim_N \frac{\displaystyle\sum_{n=1}^{\infty} \log 2^N}{\displaystyle\sum_{n=1}^{\infty} \log(2^N \pi n)} \\
&= \lim_N \frac{N \log 2 \cdot \displaystyle\sum^{\infty} 1}{(N \log 2 + \log \pi) \displaystyle\sum^{\infty} 1 + \cancel{\displaystyle\sum_{n=1}^{\infty} \log n}} \qquad (16) \\
&\sim \lim_N \frac{N \log 2}{N \log 2 + \log \pi} = 1,
\end{aligned}$$

where the logarithmic series $\sum_n \log n$ can be neglected (although it converges, in the *zeta-regularization* sense, to ½log2$\pi$) as its partial sums increase as $\mathrm{O}(\log n)$, i.e. slower than any linear series, such as $\sum 1$.

The ISS of *Cantor set* $\mathrm{C} \subset [0,1]$ (cfr. [3], [5]) is made up of $w_1(x) = \frac{1}{3}x$ and $w_2(x) = \frac{1}{3}(x+2)$. Nonzero spectrum is $\{-3^N n\pi\}_{n\in\mathrm{N}}$, with $\mathrm{mul}(3^N n\pi)=2^N$:

$$\dim_\Delta \mathrm{C} = \lim_N \frac{\sum_{n=1}^{\infty} \log 2^N}{\sum_{n=1}^{\infty} \log(3^N \pi n)}$$

$$= \lim_N \frac{N\log 2 \cdot \sum^{\infty} 1}{(N\log 3 + \log\pi)\sum^{\infty} 1 + \cancel{\sum_{n=1}^{\infty} \log n}} \tag{17}$$

$$\sim \lim_N \frac{N\log 2}{N\log 3 + \log\pi} = \frac{\log 2}{\log 3} = \dim_{\mathrm{B}} \mathrm{C}.$$

Disconnected ISSs were chosen instead of 'just-touching' ones because interconnective eigenvalues were proven not to influence the spectral dimension and such examples provide an easy, closed-form computability as well.